\documentclass[10pt,reqno,twoside]{article}
\usepackage{natbib}
\usepackage{amsfonts}
\usepackage{amssymb}
\usepackage{amsbsy}
\usepackage{amsmath,cite,graphicx}
\usepackage[top=20mm,bottom=20mm,left=15mm,right=15mm]{geometry}
\newcommand\blfootnote[1]{%
  \begingroup
  \renewcommand\thefootnote{}\footnote{#1}%
  \addtocounter{footnote}{-1}%
  \endgroup
}

\numberwithin{equation}{section}
\numberwithin{figure}{section} 

\begin{document}

\title{Optimal Perturbation Iteration Method for Bratu-Type Problems}
\author{Sinan DEN\.{I}Z$^{1}$ \ and Necdet B\.{I}LD\.{I}K$^{1}$\footnote{Corresponding author}}

\date{{\small$^{1}$ Department of Mathematics, Faculty of Art and Sciences, Celal Bayar University, 45040 Manisa, Turkey.\\ e-mails: sinan.deniz@cbu.edu.tr,n.bildik@cbu.edu.tr}}

\maketitle

\blfootnote{
 \rule{0.63cm}{0cm}\emph{Mathematics Subject Classification:} 47H14, 34L30, 80A20}

\begin{abstract}
In this paper, we introduce the new optimal perturbation iteration method based on the perturbation iteration algorithms for the approximate solutions of nonlinear differential equations of many types. The proposed method is illustrated by studying Bratu-type equations.  Our results show that only a few terms are required to obtain an approximate solution which is more accurate and efficient than many other methods in the literature.
\end{abstract}

\section{Introduction}

Many nonlinear differential equations are used in many scientific studies and most of them cannot be solved analytically using traditional methods. Therefore these problems are often handled by a broad class of analytical and numerical methods such as Adomian decomposition method \citep{adomian1988review,deniz2014comparison},  Taylor collocation method \citep{bildik2015comparison}, differential transform method \citep{bildik2006use}, homotopy perturbation method \citep{ozics2008he}, variational iteration method \citep{he2003variational}. These methods can give accurate solutions to nonlinear problems but they have also some problems about the convergence region of their series solution. These regions are generally small according to the desired solution. In order to cope with this task, researchers have recently proposed some new methods \citep{marinca2008application,liao2012optimal,idrees2010application}.  Perturbation iteration method is one of them and has been recently developed by Pakdemirli et.al. It has proven that this method is very effective for solving many nonlinear equations arising in scientific world \citep{aksoy2010new,aksoy2012new,csenol2013perturbation,timuccin2013new,khalid2015solving}.
In the presented study, we construct a new optimal perturbation iteration method which is applicable to a wide range of equations and does not require special transformations. In order to show the efficiency of the proposed method, we try to solve Bratu initial and boundary value problems  which are used in a large variety of applications, such as the fuel ignition model of the theory of thermal combustion, the thermal reaction process model, radioactive heat transfer, nanotechnology and theory of chemical reaction \citep{doha2013efficient,he2014variational,raja2014solution}.

\section{Perturbation Iteration Method}

Pakdemirli and his co-workers have modified the well-known perturbation method to construct perturbation iteration method (PIM).  PIM has been efficiently applied to some strongly nonlinear systems and yields very approximate results \citep{aksoy2012new,csenol2013perturbation}. In this section; we give basic information about perturbation iteration algorithms. They are classified with respect to the number of correction terms  ($n$) and with respect to the degrees of derivatives in the Taylor expansions($m$). Briefly, this process is represented as PIA ($n,m$).
\vglue0.2cm
\textit{\underbar{PIA (1,1) }}\underbar{}

\noindent In order to illustrate the algorithm, consider a second-order differential equation in closed form:
\begin{equation} \label{2.1}
F({y}'',{y}',y,\varepsilon )=0
\end{equation}
\noindent where  $y=y(x)$ and $\varepsilon$  is the perturbation parameter. For PIA (1,1), we take one correction term from the perturbation expansion:
\begin{equation} \label{2.2}
{y_{n+1}}={y_{n}}+\varepsilon {\left( {y_{c}} \right)}_{n}
\end{equation}
\noindent Substituting (2.2) into (2.1) and then expanding in a Taylor series gives
\begin{eqnarray} \label{2.3}
F({y_{n}}^{\prime \prime },{{y}_{n}}^{\prime },{{y}_{n}},0)+{{F}_{y}}{{({{y}_{c}})}_{n}}\varepsilon +{{F}_{{{y}'}}}{{({{y}_{c}}^{\prime })}_{n}}\varepsilon+
{{F}_{{{y}''}}}{{({{y}_{c}}^{\prime \prime })}_{n}}\varepsilon +{{F}_{\varepsilon }}\varepsilon =0
\end{eqnarray}
\noindent Rearranging the Eq. (2.3) yields a linear second order differential equation:
\begin{eqnarray} \label{2.4}
{{\left( {{y}_{c}}^{\prime \prime } \right)}_{n}}+\frac{{{F}_{{{y}'}}}}{{{F}_{{{y}''}}}}{{\left( {{y}_{c}}^{\prime } \right)}_{n}}+\frac{{{F}_{y}}}{{{F}_{{{y}''}}}}{{\left( {{y}_{c}} \right)}_{n}}=-\frac{\frac{F}{\varepsilon }+{{F}_{\varepsilon }}}{{{F}_{{{y}''}}}}
\end{eqnarray}

\noindent We can easily obtain $(y_{c})_{0}$  from Eq. (2.4) by using an initial guess  ${{y}_{0}}$. Then  first approximation  ${{y}_{1}}$ is determined by using this information.
\vglue0.2cm
\textit{\underbar{PIA (1,2) }}\underbar{}

\noindent As distinct from PIA(1,1),  we need to take  $n = 1, m = 2$ to obtain PIA (1,2). In other words, second order derivatives must be taken  into consideration:
\begin{equation} 
\begin{array}{l} F(y_{n} ^{''} ,y_{n} ^{'} ,y_{n} ,0)+F_{y} (y_{c} )_{n} \varepsilon+ F_{y'}(y_{c} ^{{'} } )_{n} \varepsilon +F_{y''}(y_{c} ^{{'} {'} } )_{n} \varepsilon  + F_{\varepsilon } \varepsilon +\frac{1}{2} \varepsilon ^{2} F_{y''y''} (y_{c} ^{{'} {'} } )_{n}^{2}
+\frac{1}{2} \varepsilon ^{2} F_{y'y'}(y_{c} ^{{'} } )_{n}^{2} +\frac{1}{2} \varepsilon ^{2} F_{yy} (y_{c} )_{n}^{2} \\
+\varepsilon ^{2} F_{y''y'} (y_{c} ^{{'} {'} } )_{n} (y_{c} ^{{'} } )_{n}+\varepsilon ^{2} F_{y'y} (y_{c} ^{{'} } )_{n} (y_{c} )_{n}
+\varepsilon ^{2} F_{y''y} (y_{c} ^{{'} {'} } )_{n} (y_{c} )_{n}+F_{\varepsilon y''} (y_{c} ^{{'} {'} } )_{n} \varepsilon ^{2}
+F_{\varepsilon y'} (y_{c} ^{{'} } )_{n} \varepsilon ^{2}+F_{\varepsilon y} (y_{c} )_{n} \varepsilon ^{2}
+\frac{1}{2} \varepsilon ^{2} F_{\varepsilon \varepsilon } =0
\end{array}
\end{equation}
\noindent or by rearranging
\begin{equation} \label{2.6}
\begin{array}{l} {(y_{c} ^{''} )_{n} \left(\varepsilon F_{y''} +\varepsilon ^{2} F_{\varepsilon y''} \right)+(y_{c} ^{'} )_{n} \left(\varepsilon F_{y'} +\varepsilon ^{2} F_{\varepsilon y'} \right)+(y_{c} )_{n} \left(\varepsilon F_{y} +\varepsilon ^{2} F_{\varepsilon y} \right)+}
(y_{c} ^{''} )_{n}^{2} \left(\frac{\varepsilon ^{2}}{2} F_{y''y''} \right)
+(y_{c} ^{{'} } )_{n}^{2} \left(\frac{\varepsilon ^{2} }{2} F_{y'y'} \right)+\\(y_{c} )_{n}^{2} \left(\frac{\varepsilon ^{2} }{2} F_{yy} \right)+
(y_{c} )_{n} (y_{c} ^{'})_{n} \left(\varepsilon ^{2} F_{y'y} \right)+
(y_{c} ^{{'} {'} } )_{n} (y_{c} ^{{'} } )_{n} \left(\varepsilon ^{2} F_{y'y''} \right)+
(y_{c} ^{{'} {'} } )_{n} (y_{c} )_{n} \left(\varepsilon ^{2} F_{yy''} \right)=-F-F_{\varepsilon } \varepsilon -\frac{\varepsilon ^{2} F_{\varepsilon \varepsilon } }{2}.  \end{array}
\end{equation}
\noindent Note that all derivatives and functions  are calculated at $\varepsilon=0$ . By means of (2.2) and (2.6), iterative scheme is developed for the equation under consideration.

\section{Optimal Perturbation Iteration Method}

To illustrate the basic concept of the optimal perturbation iteration method (OPIM), we first reconsider the Eq. (2.1) as:

\begin{equation} \label{3.1}
F({y}'',{y}',y,\varepsilon )=Ly+N({y}'',{y}',y,\varepsilon ), ~~~~ B(y,y')=0
\end{equation}

\noindent where $L$  is a linear operator, $N$ denotes the nonlinear terms and $B$ is a boundary operator respectively. We then expand only nonlinear terms in a Taylor series to decrease the volume of calculations. Because, it is useless and unnecessary to expand the whole equation for each problems.  This is the first step of OPIM to decrease the time needed for computations. \\

\noindent After the Eqs. (2.4) and (2.6) in the solution processes for PIAs (1, m), we offer to use the formula
\begin{equation} \label{3.2}
{{y}_{n+1}}={{y}_{n}}+{{S}_{n}}(\varepsilon ){{\left( {{y}_{c}} \right)}_{n}}
\end{equation}

\noindent to increase the accuracy of the results and effectiveness of the method. Here $S_{n}(\varepsilon)$   is an auxiliary function which provides us to adjust and control the convergence. This is the crucial point of OPIM. The choices of functions $S_{n}(\varepsilon)$   could be exponential, polynomial, etc. In this study, we select auxiliary function in the form
\begin{equation} \label{3.3}
{{S}_{n}}(\varepsilon )={{C}_{0}}+\varepsilon {{C}_{1}}+{{\varepsilon }^{2}}{{C}_{2}}+{{\varepsilon }^{3}}{{C}_{3}}+\cdots =\sum\limits_{i=0}^{n}{{{\varepsilon }^{i}}{{C}_{i}}}
\end{equation}

\noindent where  $C_0,C_1,\ldots $  are constants which are to be determined later. \\
\noindent  The following algorithm can be used for OPIM:

\textbf{a)} Take the governing differential equation as:
\begin{equation} \label{3.4}
Ly+N({y}'',{y}',y,\varepsilon )=0,y=y(x),a\le x\le b	
\end{equation}

\textbf{b)} Substitute (2.2) into the nonlinear part of (3.4) and expand it in a Taylor series:

\begin{equation} \label{3.5}
N({{y}_{n}}^{\prime \prime },{{y}_{n}}^{\prime },{{y}_{n}},0)+{{N}_{y}}{{({{y}_{c}})}_{n}}\varepsilon +{{N}_{{{y}'}}}{{({{y}_{c}}^{\prime })}_{n}}\varepsilon +{{N}_{{{y}''}}}{{({{y}_{c}}^{\prime \prime })}_{n}}\varepsilon +{{N}_{\varepsilon }}\varepsilon =0
\end{equation}
\noindent and
\begin{equation} \label{3.6}
\begin{array}{l} N+{{N}_{y}}{{({{y}_{c}})}_{n}}\varepsilon +{{N}_{{{y}'}}}{{({{{y}'}_{c}})}_{n}}\varepsilon +{{N}_{\varepsilon }}\varepsilon +{{N}_{\varepsilon y}}{{({{y}_{c}})}_{n}}{{\varepsilon }^{2}}+{{N}_{\varepsilon {y}'}}{{({{{y}'}_{c}})}_{n}}{{\varepsilon }^{2}}+
\frac{{{N}_{\varepsilon \varepsilon }}{{\varepsilon }^{2}}}{2}+\frac{{{N}_{yy}}{{\varepsilon }^{2}}({{y}_{c}})_{n}^{2}}{2}+\frac{{{N}_{{y}'{y}'}}{{\varepsilon }^{2}}({{{{y}'}}_{c}})_{n}^{2}}{2}=0 \end{array}
\end{equation}

\textbf{c)} After finding $(y_c)_{0}$  for each algorithms as in PIAs(1,m), substitute it into Eq.(3.2) to find the first approximate result:
\begin{equation} \label{3.7}
{{y}_{1}}={{y}_{0}}+{{S}_{0}}(\varepsilon ){{\left( {{y}_{c}} \right)}_{0}}={{y}_{0}}+{{C}_{0}}{{\left( {{y}_{c}} \right)}_{0}}
\end{equation}

\noindent By using initial condition and setting $\varepsilon=1$   yields

\begin{equation} \label{3.8}
{{y}_{1}}=y(x,{{C}_{0}})
\end{equation}

\noindent Using the Eq.(3.8) and repeating the similar steps, we have:

\begin{equation} \label{3.9}
\begin{array}{l} {y_2(x,{{C}_{0}},{{C}_{1}})={{y}_{1}}+{{S}_{1}}(\varepsilon ){{\left( {{y}_{c}} \right)}_{0}}={{y}_{1}}+\left( {{C}_{0}}+{{C}_{1}} \right){{\left( {{y}_{c}} \right)}_{1}}} \\

 y_3(x,{{C}_{0}},{{C}_{1}},{{C}_{2}})={{y}_{2}}+\left( {{C}_{0}}+{{C}_{1}}+{{C}_{2}} \right){{\left( {{y}_{c}} \right)}_{2}} \\

 \qquad \vdots \\

 y_m(x,{{C}_{0}},\ldots ,{{C}_{m-1}})={{y}_{m-1}}+\left( {{C}_{0}}+\cdots +{{C}_{m-1}} \right){{\left( {{y}_{c}} \right)}_{m-1}}   \end{array}
\end{equation}

\textbf{d)} Substitute the approximate solution $y_m$ into the Eq.(3.4)  and the general problem results in the following residual:
\begin{equation} \label{3.10}
R(x,{{C}_{0}},\ldots ,{{C}_{m-1}})=L\left( {{y}_{m}}(x,{{C}_{0}},\ldots ,{{C}_{m-1}}) \right)+N\left( {{y}_{m}}(x,{{C}_{0}},\ldots ,{{C}_{m-1}}) \right)
\end{equation}
\noindent Obviously, when $R(x,{{C}_{0}},\ldots ,{{C}_{m-1}})=0$ then the approximation ${{y}_{m}}(x,{{C}_{0}},\ldots ,{{C}_{m-1}})={{y}^{(m)}}(x,{{C}_{i}})$ will be the exact solution. Generally it doesn't happen, especially in nonlinear equations. To determine the optimum values of ${{C}_{0}},{{C}_{1}},\ldots $; we here use the equations

\begin{equation} \label{3.11}
R({{x}_{1}},{{C}_{i}})=R({{x}_{2}},{{C}_{i}})=\cdots =R({{x}_{m}},{{C}_{i}})=0,i=0,1,\ldots ,m-1
\end{equation}

\noindent where ${{x}_{i}}\in (a,b)$.  Generally it is quite impossible to solve the system of Eqs. (3.11) other than numerically. Therefore, one needs to use a computer program such that Mathematica, Maple etc. Note that the solution of the system (3.11) is not unique,  but all obtained constants would yield the same approximate solutions.

The constants ${{C}_{0}},{{C}_{1}},\ldots $  can also be defined from the method of least squares:
\begin{equation} \label{3.12}
J({{C}_{0}},\ldots ,{{C}_{m-1}})=\int\limits_{a}^{b}{{{R}^{2}}}(x,{{C}_{0}},\ldots ,{{C}_{m-1}})dx
\end{equation}

\noindent where $a$ and $b$  are selected from the domain of the problem. Putting these constants into the last one of the Eqs. (3.9), the approximate solution of order $m$ is well-determined. It should be also emphasized that, the Eq. (3.12) is not always useful to find the constants ${{C}_{0}},{{C}_{1}},\ldots $ especially for strongly nonlinear equations. So,we use the Eq. (3.11) to get those constants in this work. For much more information and different usage about this process, please see \citep{herisanu2015analytical,marinca2012optimal}

\section{Applications}

\textbf{Example 1.} Consider the following nonlinear differential equation \citep{wazwaz2005adomian}:
\begin{equation} \label{4.1}
{y}''-2{{e}^{y}}=0,y(0)={y}'(0)=0,0\le x\le 1.
\end{equation}

\noindent which has the exact solution $y=-2\ln \left( \cos x \right)$.

\textit{\underbar{OPIA (1,1) }}\underbar{}

\noindent Consider the Eq. (4.1) as:
\begin{equation} \label{4.2}
F({y}'',y,\varepsilon )={y}''-2{{e}^{\varepsilon y}}={y}''+N(y,\varepsilon ).
\end{equation}

\noindent  OPIA (1,1)  requires to compute:
\begin{equation} \label{4.3}
N({{y}_{n}},0)+{{N}_{y}}({{y}_{n}},0){{({{y}_{c}})}_{n}}\varepsilon +{{N}_{\varepsilon }}\varepsilon =0
\end{equation}

\noindent which is approximately half of the volume of calculations that in PIA(1,1).
\noindent Using the Eqs. (2.2), (4.3) and setting  $\varepsilon=1$ yields

\begin{equation} \label{4.4}
{{({{y}_{c}}^{\prime \prime })}_{n}}=-{{y}_{n}}^{\prime \prime }+2{{y}_{n}}+2
\end{equation}

\noindent One may start the iteration by taking a trivial solution which satisfies the given initial conditions:
\begin{equation} \label{4.5}
{y}_{0}=0.
\end{equation}

\noindent Substituting (4.5) into the Eq. (4.4), we have

\begin{equation} \label{4.6}
{{\left( {{y}_{c}} \right)}_{0}}={{x}^{2}}+c
\end{equation}

\noindent Now, Eq. (4.6) is inserted into Eq. (3.2) and applying the initial conditions we get
\begin{equation}\label{4.7}
{{y}_{1}}={{y}_{0}}+{{S}_{n}}(\varepsilon ){{\left( {{y}_{c}} \right)}_{0}}={{C}_{0}}{{x}^{2}}
\end{equation}

\noindent It is worth mentioning that $y_1$  does not represent the first correction term; rather it is the approximate solution after the first iteration. Following the same procedure, we obtain new and more approximate results:

\begin{equation}\label{4.8}
{{y}_{2}}={{C}_{0}}{{x}^{2}}+({{C}_{0}}+{{C}_{1}})({{x}^{2}}-{{C}_{0}}{{x}^{2}}+\frac{{{C}_{0}}{{x}^{4}}}{6})
\end{equation}

\begin{equation}\label{4.9}
\begin{array}{l}
   {{{y}_{3}}=\left[ 2{{C}_{0}}+{{C}_{1}}-{{C}_{0}}({{C}_{0}}+{{C}_{1}})+(-1+{{C}_{0}})(-1+{{C}_{0}}+{{C}_{1}})({{C}_{0}}+{{C}_{1}}+{{C}_{2}}) \right]{{x}^{2}}} \\
  {+\left[ \frac{1}{6}{{C}_{0}}({{C}_{0}}+{{C}_{1}})+\frac{1}{15}\left( -5{{C}_{0}}^{2}-5{{C}_{0}}(-1+{{C}_{1}})+\frac{5{{C}_{1}}}{2} \right)({{C}_{0}}+{{C}_{1}}+{{C}_{2}}) \right]{{x}^{4}}}
  {+\left[ \frac{1}{90}{{C}_{0}}({{C}_{0}}+{{C}_{1}})({{C}_{0}}+{{C}_{1}}+{{C}_{2}}) \right]{{x}^{6}}}
\end{array}
\end{equation}

\noindent To determine the constants, we proceed as in section 3. First, the residual
\begin{equation}\label{4.10}
\begin{array}{l}
   R(x,{{C}_{0}},{{C}_{1}},{{C}_{2}})=L\left( {{y}_{3}}(x,{{C}_{0}},{{C}_{1}},{{C}_{2}}) \right)+N\left( {{y}_{3}}(x,{{C}_{0}},{{C}_{1}},{{C}_{2}}) \right)=

   2{{C}_{0}}+({{C}_{0}}+{{C}_{1}})(2-2{{C}_{0}}+2{{C}_{0}}{{x}^{2}}) +\frac{({{C}_{0}}+{{C}_{1}}+{{C}_{2}})}{15}\times \\ \left[ 30(-1+{{C}_{0}})(-1+{{C}_{0}}+{{C}_{1}})+12\left( -5{{C}_{0}}^{2}-5{{C}_{0}}(-1+{{C}_{1}})+\frac{5{{C}_{1}}}{2} \right){{x}^{2}}+5{{C}_{0}}({{C}_{0}}+{{C}_{1}}){{x}^{4}} \right] \\
  -2Exp\left[ \begin{array}{l}
   {{C}_{0}}{{x}^{2}}+({{C}_{0}}+{{C}_{1}})\left( {{x}^{2}}-{{C}_{0}}{{x}^{2}}+\frac{{{C}_{0}}{{x}^{4}}}{6} \right)+ +\frac{1}{6}{{C}_{0}}({{C}_{0}}+{{C}_{1}}){{x}^{6}}\\
  \frac{({{C}_{0}}+{{C}_{1}}+{{C}_{2}})}{15}\left( \begin{array}{l}
   15(-1+{{C}_{0}})(-1+{{C}_{0}}+{{C}_{1}}){{x}^{2}}
  +\left( -5{{C}_{0}}^{2}-5{{C}_{0}}(-1+{{C}_{1}})+\frac{5{{C}_{1}}}{2} \right){{x}^{4}} \\
\end{array} \right) \\
\end{array} \right] \\
\end{array}
\end{equation}

\noindent is constructed for the third order approximation. Using the Eq. (3.11)  with $x=0.3,0.6,0.9$ , we get
\begin{equation}\label{4.11}
{{C}_{0}}=1.00096007239,{{C}_{1}}=0.034138423506,{{C}_{2}}=-0.049127633506
\end{equation}

\noindent Inserting the constants into the Eq. (4.9), we obtain the approximate solution of the third order:

\begin{equation}\label{4.12}
{{y}_{3}}(x)=1.00112456947{{x}^{2}}+0.152984774463{{x}^{4}}+0.076778117636{{x}^{6}}
\end{equation}

\noindent Note that some complex numbers arise from solving the Eq. (4.10). They can also be used instead of  $C_0,C_1,C_2$ to get the same result. We here give only real solutions for simplicity.

\textit{\underbar{OPIA (1,2) }}\underbar{}

\noindent One can construct the OPIA(1,2) by taking one correction term in the perturbation expansion and two derivatives in the Taylor series. Note that one needs to enter the data in Eq. (2.4) into the computer for PIA(1,2). But, it is sufficient to use

\begin{equation}\label{4.13}
N+{{N}_{y}}{{({{y}_{c}})}_{n}}\varepsilon +{{N}_{\varepsilon }}\varepsilon +{{N}_{\varepsilon y}}{{({{y}_{c}})}_{n}}{{\varepsilon }^{2}}+\frac{{{N}_{\varepsilon \varepsilon }}{{\varepsilon }^{2}}}{2}+\frac{{{N}_{yy}}{{\varepsilon }^{2}}({{y}_{c}})_{n}^{2}}{2}=0
\end{equation}

\noindent for OPIA(1,2). After making the relevant calculations, the algorithm takes the simplified form:
\begin{equation}\label{4.14}
{{({{y}_{c}}^{\prime \prime })}_{n}}-2{{({{y}_{c}})}_{n}}=-{{y}_{n}}^{\prime \prime }+2{{y}_{n}}+{{y}_{n}}^{2}+2
\end{equation}
\noindent Using the trivial solution ${{y}_{0}}=0$ , we have second order problem
\begin{equation}\label{4.15}
{{({{y}_{c}}^{\prime \prime })}_{0}}-2{{({{y}_{c}})}_{0}}=2
\end{equation}
\noindent Using the Eqs. (3.2), (4.15) and the initial conditions, we obtain
\begin{equation}\label{4.16}
{{y}_{1}}={{C}_{0}}\left( \cosh {\left( \sqrt{2}x \right)}-1 \right)
\end{equation}

\noindent Following the same procedure using (4.16), the second iteration is obtained as

\begin{equation}\label{4.17}
{{y}_{2}}=\frac{1}{3}\left( \begin{array}{l}
   -3{{C}_{0}}+3{{C}_{0}}\cosh[\sqrt{2}x]-\frac{3{{C}_{0}}^{2}({{C}_{0}}+{{C}_{1}})x\sinh[\sqrt{2}x]}{\sqrt{2}} \\
  +\left( ({{C}_{0}}+{{C}_{1}})(6+{{C}_{0}}(-6+5{{C}_{1}})+{{C}_{0}}^{2}\cosh[\sqrt{2}x]) \right)\sinh{{[\frac{x}{\sqrt{2}}]}^{2}} \\
\end{array} \right)
\end{equation}

\noindent One can easily realize that, we have functional expansion for OPIA (1,2)  instead of a polynomial expansion.\\
Following the same procedure, from the residual

\begin{equation}\label{4.18}
\begin{array}{l}
   R(x,{{C}_{0}},{{C}_{1}})=L(y_2)+N(y_2)=
 \frac{1}{3}\left[ \begin{array}{l}
  -2\left( -3{{C}_{1}}+{{C}_{0}}\left( 3\left( -2+{{C}_{1}} \right)+{{C}_{0}}\left( 3+{{C}_{0}}+{{C}_{1}} \right) \right) \right)\cosh[\sqrt{2}x] \\
 +{{C}_{0}}^{2}({{C}_{0}}+{{C}_{1}})\left( 2\cosh[2\sqrt{2}x]-3\sqrt{2}x\sinh[\sqrt{2}x] \right) \\
\end{array} \right] \\
  -2Exp\left[ \frac{1}{3}\left( \begin{array}{l}
  ({{C}_{0}}+{{C}_{1}})\left( 6+{{C}_{0}}(-6+5{{C}_{0}})+{{C}_{0}}^{2} \cosh[\sqrt{2}x] \right)sinh{{[\frac{x}{\sqrt{2}}]}^{2}} \\
 -\frac{3{{C}_{0}}^{2}({{C}_{0}}+{{C}_{1}})x\sinh[\sqrt{2}x]}{\sqrt{2}}  -3{{C}_{0}}+3{{C}_{0}}\cosh[\sqrt{2}x]
\end{array} \right) \right]
\end{array}
\end{equation}

\noindent the constants $C_0$ and $C_1$   can be determined as
\begin{equation}\label{4.19}
{{C}_{0}}=1.000861120478,{{C}_{1}}=0.0266135748038
\end{equation}

 \noindent Thus, we have the second-order approximate solution:

\begin{equation}\label{EQ38}
\begin{array}{l}{
  {{y}_{2}}(x)=-1.1784311655118591x\sinh (\sqrt{2}x)+2.1406095945289634\cosh (\sqrt{2}x)} \\
  {+0.13215241067298802\cosh (2\sqrt{2}x)-2.2728821098149927 }
\end{array}
\end{equation}

\noindent One can also compute more approximate results by following the same procedure with a computer program. We do not give higher iterations due to huge amount of calculations. Figure 1 and Table 1 show a comparison of OPIAs and exact solution. It is clear that the results obtained by OPIM are more accurate than those of PIM in \citep{aksoy2010new}.

\begin{table}\caption{Comparison of absolute errors of  Example 1 at different orders of approximations.}
\tabcolsep 5.8pt
\centering
\begin{tabular}{|p{0.15in}|p{0.8in}|p{0.8in}|p{0.8in}|p{0.8in}|p{1in}|p{0.8in}|} \hline
       \textit{x} & \multicolumn{3}{|p{3in}|}{Absolute  errors for OPIA(1,1)   solutions } & \multicolumn{2}{|p{2in}|}{Absolute  errors for  OPIA(1,2)   solutions } & Exact solution \\ \hline
 &    $\left|y-y_{1} \right|$ &    $\left|y-y_{2} \right|$ &      $\left|y-y_{3} \right|$ &      $\left|y-y_{1} \right|$ &     $\left|y-y_{2} \right|$ & $y=-2\ln \left(\cos x\right)$ \\ \hline
0.1 & 0.000449452 & 0.000169553 & 9.9097$\times$10${}^{-6}$ & 2.402$\times$10${}^{-6}$ & 9.4728$\times$10${}^{-6}$ & $0.010016711$ \\
0.2 & 0.001595127 & 0.000583911 & 2.5126$\times$10${}^{-5}$ & 9.453$\times$10${}^{-6}$ & 3.3152$\times$10${}^{-5}$ & $0.040269546$ \\
0.3 & 0.002812140 & 0.000976872 & 1.3047$\times$10${}^{-5}$ & 1.9420$\times$10${}^{-5}$ & 2.7254$\times$10${}^{-5}$ & $0.091383311$ \\
0.4 & 0.003000543 & 0.000963454 & 4.7213$\times$10${}^{-5}$ & 2.4899$\times$10${}^{-6}$ & 4.4563$\times$10${}^{-6}$ & $0.164458038$ \\
0.5 & 0.000485555 & 0.000139394 & 0.000126132 & 4.916$\times$10${}^{-6}$ & 5.55112$\times$10${}^{-8}$ & $0.261168480$ \\
0.6 & 0.007148548 & 0.001752633 & 0.000116507 & 8.8755$\times$10${}^{-5}$ & 7.2047$\times$10${}^{-5}$ & $0.383930338$ \\
0.7 & 0.023329621 & 0.004551758 & 0.000144037 & 0.000354849 & 7.0044$\times$10${}^{-5}$ & $0.536171515$ \\
0.8 & 0.052947212 & 0.007229526 & 0.000727717 & 0.000982654 & 0.000128213 & $0.722781493$ \\
0.9 & 0.103126097 & 0.007116353 & 0.001202366 & 0.002323371 & 0.000452361 & $0.950884887$ \\
  1 & 0.184637089 & 0.001509956 & 0.000365479 & 0.005024005 & 4.44089$\times$10${}^{-8}$ & $1.231252940$ \\ \hline
\end{tabular}
\end{table}

\begin{figure}[h]
\centerline{\includegraphics[width=0.7\textwidth]{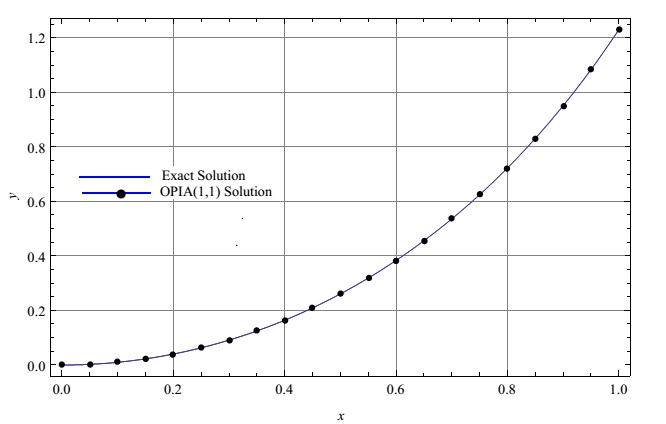}}
\caption{Comparison between the three-term OPIA(1,1) approximate solution and the exact solution for Example 1.}
\end{figure}

\vglue0.3cm
\noindent \textbf{Example 2.} Bratu's  first  boundary value problem is given as \citep{wazwaz2005adomian}:
\begin{equation} \label{EQ39}
y''+\lambda e^{y} {\rm =0,\; \; \; 0}\le x\le 1,{\rm \; \; \; \; }y(0)=y(1)=0
\end{equation}
\noindent with the exact solution $y(x)=-2\ln \left[\frac{\cosh \left(\left(x-\frac{1}{2} \right)\frac{\theta }{2} \right)}{\cosh \left(\frac{\theta }{4} \right)} \right]$

\noindent where $\theta $ satisfies
\[\theta =\sqrt{2\lambda } \cosh \left(\frac{\theta }{4} \right).\]
\textit{\underbar{OPIA (1,1) }}

\noindent \textit{\underbar{}}

\noindent An artificial perturbation parameter is inserted for the Eq. (4.21) as follows
\begin{equation} \label{EQ40}
F(y'',y,\varepsilon )=y''+\lambda e^{\varepsilon y} =Ly+N(y,\varepsilon )=0.
\end{equation}
\noindent By making necessary computations using the Eqs. (2.2), (4.3) and  setting $\varepsilon =1$, we easily get
\begin{equation} \label{GrindEQ__41_}
(y_{c} ^{{'} {'} } )_{n} =-\lambda y_{n} -\left(y_{n} ^{{'} {'} } +\lambda \right).
\end{equation}
\noindent One may start with the trivial solution
\begin{equation} \label{EQ42}
y_{0} =0
\end{equation}

\noindent and using the Eq.(3.2)  the iterations are reached as follows:
\begin{equation} \label{EQ43}
y_{1} =-\frac{\lambda C_{0} }{2} \left(x^{2} -x\right)
\end{equation}

 \begin{equation}\label{EQ44}
  y_{2} =-\frac{\lambda C_{0} }{2} (x^{2} -x)+\frac{(C_{0} +C_{1} )}{24} (-1+x)x\lambda \left[-12+\left(12+(-1-x+x^{2} )\lambda \right)C_{0} \right]
 \end{equation}

 \begin{equation} \label{EQ45}
\begin{array}{l} {y_{3} =  \frac{\lambda x}{720}\times}\\ {{\left[\begin{array}{l} {30(C_{0} +C_{1} )x\left(-12+C_{0} (12+\lambda (-2+x)x)\right)-360C_{0} (-1+x)-} \\ {(C_{0} +C_{1} +C_{2} )x\left\{\begin{array}{l} {360(-1+C_{0} )(-1+C_{0} +C_{1} )-60C_{0} (-1+C_{0} +C_{1} )\lambda x} \\ {+30\left(-C_{1} +2C_{0} (-1+C_{0} +C_{1} )\right)\lambda x^{2} -3C_{0} (C_{0} +C_{1} )\lambda ^{2} x^{3}} \\ {+C_{0} (C_{0} +C_{1} )\lambda ^{2} x^{4} } \end{array}\right\}} \\ {{\rm \; \; \; \; \; \; \; \; \; \; \; \; \; \; \; \; \; \; \; \; \; \; \; \; \; \; \; }} \end{array}\right]}} \end{array}
\end{equation}

\noindent For the constants $C_{0} ,C_{1} $  and $C_{2} $ , the method given in Section 3.(d) is used, and we obtain the following values for $\lambda =1$:
\begin{equation}\label{EQ46}
C_{0} =0.00896621251,C_{1} =0.086955412771,C_{2} =-0.000213669444
\end{equation}

\noindent for the $x_i=0.3,0.6,0.9$. Thus, the approximate solution of the third order is:
\begin{equation} \label{EQ47}
\begin{array}{l} {y_{3} (x)=0.549359811237294x-0.5001682773565349x^{2} -0.09044598683883211x^{3} } \\ {+0.025373292932254158x^{4} +0.023828818721721115x^{5} -0.00794805048551452x^{6} } \end{array}
\end{equation}

\noindent

\noindent \textit{\underbar{OPIA (1,2) :}}

\noindent One just needs to construct
\begin{equation} \label{EQ48}
N+N_{y} (y_{c} )_{n} \varepsilon +N_{\varepsilon } \varepsilon +N_{\varepsilon y} (y_{c} )_{n} \varepsilon ^{2} +\frac{N_{\varepsilon \varepsilon } \varepsilon ^{2} }{2} +\frac{N_{yy} \varepsilon ^{2} (y_{c} )_{n}^{2} }{2} =0
\end{equation}
\noindent where
\begin{equation} \label{EQ49}
N(y,\varepsilon )=\lambda e^{\varepsilon y}.
\end{equation}

\noindent After making the relevant calculations, the algorithm  takes the simplified form:
\begin{equation} \label{EQ50}
(y_{c} ^{{'} {'} } )_{n} +\lambda (y_{c} )_{n} =-\lambda y_{n} -y_{n} ^{{'} {'} } -\lambda -\frac{\lambda }{2} y_{n}^{2}  .
\end{equation}

\noindent Using the Eqs. (3.2), (4.24),(4.32) and the initial conditions, we obtain
\begin{equation} \label{EQ51}
y_{1} =C_{0} \left(-1+ \cos[x\sqrt{\lambda } ]+ \sin[x\sqrt{\lambda } ] \tan[\frac{\sqrt{\lambda } }{2} ]\right)
\end{equation}

\begin{equation} \label{EQ52}
\begin{array}{l} {y_{2} =C_{0} \left(-1+\cos[x\sqrt{\lambda } ]+\sin[x\sqrt{\lambda } ]\tan[\frac{\sqrt{\lambda } }{2} ]\right)+\frac{(C_{0} +C_{1} )}{48} \left(\sec[\frac{\sqrt{\lambda } }{2} ]^{2} (-24(1-C_{0} +C_{0} ^{2} )\cos[x\sqrt{\lambda } ]^{2} \right)} \\

 {+(C_{0} +C_{1} )\left(\begin{array}{l} {(12-12C_{0} +C_{0} ^{2} )\cos[\sqrt{\lambda } ]-C_{0} ^{2} \cos[(1-3x)\sqrt{\lambda } ]+3(4-4C_{0} +3C_{0} ^{2} } \\ {+C_{0} ^{2} \cos[(1-2x)\sqrt{\lambda } ]+2(-1+C_{0} )\cos[(-1+x)\sqrt{\lambda } ]+C_{0} ^{2} \cos[2x\sqrt{\lambda } ]} \\ {-2\cos[(1+x)\sqrt{\lambda } ]+2C_{0} \cos[(1+x)\sqrt{\lambda } ]-C_{0} ^{2} \cos[(1+x)\sqrt{\lambda } ]} \\ {-2C_{0} ^{2} x\sqrt{\lambda } \sin[\sqrt{\lambda } ])} \end{array}\right)} \\ {+2(C_{0} +C_{1} )\cos[x\sqrt{\lambda } ]+(C_{0} +C_{1} )\sec[\frac{\sqrt{\lambda } }{2} ] \sin[x\sqrt{\lambda } ]\times } \\ {\left[\begin{array}{l} {6C_{0} ^{2} (-2+3x)\sqrt{\lambda } \cos[\frac{\sqrt{\lambda } }{2} ]+6C_{0} ^{2} x\sqrt{\lambda } \cos[\frac{3\sqrt{\lambda } }{2} ]+12\sin[\frac{\sqrt{\lambda } }{2} ]-12C_{0} \sin[\frac{\sqrt{\lambda } }{2} ]+17C_{0} ^{2} \sin[\frac{\sqrt{\lambda } }{2} ]} \\
{+12\sin[\frac{3\sqrt{\lambda } }{2} ]-12C_{0} \sin[\frac{3\sqrt{\lambda } }{2} ]+C_{0} ^{2} \sin[\frac{3\sqrt{\lambda } }{2} ]+C_{0} ^{2} Sin[\frac{1}{2} (1-6x)\sqrt{\lambda } ]-6C_{0} ^{2} \sin[\frac{1}{2} (1-4x)\sqrt{\lambda } ]} \\
{-3C_{0} ^{2} \sin[\frac{1}{2} (3-4x)\sqrt{\lambda } ]+18\sin[\frac{1}{2} (1-2x)\sqrt{\lambda } ]-18C_{0} \sin[\frac{1}{2} (1-2x)\sqrt{\lambda } ]+18C_{0} ^{2} \sin[\frac{1}{2} (1-2x)\sqrt{\lambda } ]} \\ {+C_{0} ^{2} \sin[\frac{3}{2} (1-2x)\sqrt{\lambda } ]+6\sin[\frac{1}{2} (3-2x)\sqrt{\lambda } ]-6C_{0} \sin[\frac{1}{2} (3-2x)\sqrt{\lambda } ]+6C_{0} ^{2} \sin[\frac{1}{2} (3-2x)\sqrt{\lambda } ]-} \\ {18\sin[\frac{1}{2} (1+2x)\sqrt{\lambda } ]+18C_{0} \sin[\frac{1}{2} (1+2x)\sqrt{\lambda } ]-15C_{0} ^{2} \sin[\frac{1}{2} (1+2x)\sqrt{\lambda } ]} \\ {-6\sin[\frac{1}{2} (3+2x)\sqrt{\lambda } ]+6C_{0} \sin[\frac{1}{2} (3+2x)\sqrt{\lambda } ]-3A^{2} \sin[\frac{1}{2} (3+2x)\sqrt{\lambda } ]+3C_{0} ^{2} \sin[\frac{1}{2} (1+4x)\sqrt{\lambda } ]} \end{array}\right]} \end{array}
\end{equation}

\noindent For the constants $C_{0} $ and $C_{1} $  in Eq.(3.34), we proceed as earlier and get
\begin{equation} \label{EQ53}
C_{0} =-1.0002036577189,C_{1} =0.099502786321
\end{equation}
for $\lambda =1$.Thus, we have the second-order approximate solution:
\begin{equation}\label{EQ54}
\begin{array}{l}
  {{y}_{2}}(x)=-1.078485122090-0.004293531433x+1.105765327206\cos[x]
 -0.0279349653844\cos[2x]+\\
 0.00065493410502\cos[3x]
 +0.6114581430450\sin[x]-0.091877388999\cos[x]\sin[x]
 +0.011352766729676\sin[3x]
\end{array}
\end{equation}

\noindent for OPIA(1,2). It can be readily seen from Figure 2 and Table 2 that approximate solutions obtained by the OPIAs are identical with that given by the analytical methods \citep{wazwaz2005adomian}.  Note that more components in the solution series can be computed to enhance the approximation.

\begin{table}\caption{Comparison of absolute errors of  Example 2 at different orders of approximations.}
\tabcolsep 5.8pt
\begin{tabular}{|p{0.15in}|p{1in}|p{1in}|p{1in}|p{1in}|p{1in}|p{0.8in}|} \hline
       \textit{x} & \multicolumn{3}{|p{3.3in}|}{ Absolute  errors for OPIA(1,1)   solutions} &
       \multicolumn{2}{|p{2in}|}{Absolute  errors for  OPIA(1,2)   solutions} & Exact solution \\ \hline

 &    $\left|y-y_{1} \right|$ & $\left|y-y_{2} \right|$ &      $\left|y-y_{3} \right|$ &      $\left|y-y_{1} \right|$ &     $\left|y-y_{2} \right|$ &    for $\lambda =1$ \\ \hline
0.1 & 1.05236$\times$10${}^{-6}$ & 8.05698$\times$10${}^{-7}$ & 1.19748$\times$10${}^{-7}$ & 5.22201$\times$10${}^{-10}$ & 1.23154$\times$10${}^{-16}$ & 0.0498465 \\
0.2 & 1.08547$\times$10${}^{-5}$ & 7.5067 $\times$ 10${}^{-7}$ & 3.35942$\times$10${}^{-8}$ & 7.90215$\times$10${}^{-9}$ & 2.36014$\times$10${}^{-15}$ & 0.0891894 \\
0.3 & 4.96318$\times$10${}^{-5}$ & 1.00521$\times$10${}^{-6}$ & 1.12813$\times$10${}^{-8}$ & 5.20476$\times$10${}^{-9}$ & 5.10365$\times$10${}^{-13}$ &  0.1176084 \\
0.4 & 9.55681$\times$10${}^{-5}$ & 5.96014$\times$10${}^{-8}$ & 9.08115$\times$10${}^{-9}$ & 2.63391$\times$10${}^{-11}$ & 5.30158$\times$10${}^{-15}$ & 0.1347894 \\
0.5 & 7.56419$\times$10${}^{-6}$ & 7.22085$\times$10${}^{-8}$ & 7.33394$\times$10${}^{-10}$ & 9.89661$\times$10${}^{-10}$ & 5.60972$\times$10${}^{-15}$& 0.1405383 \\
0.6 & 0.000121368 & 5.00123$\times$10${}^{-7}$ & 1.13418$\times$10${}^{-9}$ & 2.00569$\times$10${}^{-11}$ & 9.12054$\times$10${}^{-13}$ & 0.1347894 \\
0.7 & 0.000802364 & 4.20161$\times$10${}^{-6}$ & 6.13948$\times$10${}^{-9}$ & 4.11057$\times$10${}^{-11}$ & 2.03606$\times$10${}^{-13}$ & 0.1176084 \\
0.8 & 0.000110879 & 1.00907$\times$10${}^{-5}$ & 1.00907$\times$10${}^{-8}$ & 8.05698$\times$10${}^{-10}$ & 7.45236$\times$10${}^{-12}$ & 0.0891894 \\
0.9 & 0.000569203 & 2.10102$\times$10${}^{-5}$ & 7.75262$\times$10${}^{-8}$ & 2.05471$\times$10${}^{-9}$ & 1.00612$\times$10${}^{-12}$ &  0.0498465 \\
  \hline
\end{tabular}
\end{table}

\begin{figure}[h]
\centerline{\includegraphics[width=0.7\textwidth]{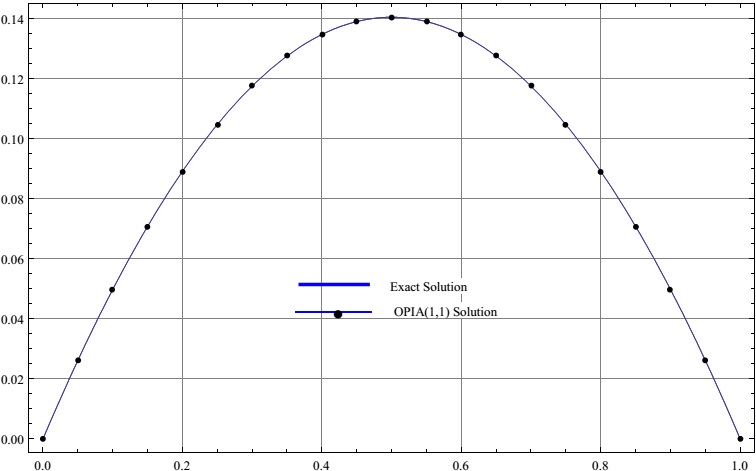}}
\caption{Comparison between the three-term OPIA(1,1) approximate solution and the exact solution for Example 2.}
\end{figure}

\vglue0.3cm
\noindent \textbf{Example 3. } Consider Bratu's second  boundary value problem \citep{wazwaz2005adomian}
\begin{equation} \label{EQ55}
y''+\pi ^{2} e^{-y} {\rm =0,\; \; \; 0}\le x\le 1,{\rm \; \; \; \; }y(0)=y(1)=0.
\end{equation}
\noindent Exact solution of this problem is mistakenly given as
\begin{equation} \label{EQ56}
y(x)=\ln \left[1+\sin \left(1+\pi x\right)\right]
\end{equation}
\noindent in  \citep{wazwaz2005adomian,batiha2010numerical}, whereas the correct exact solution is
\begin{equation} \label{EQ57}
y(x)=\ln \left[1+\sin \left(\pi x\right)\right].
\end{equation}

\textit{\underbar{OPIA (1,1) }}\underbar{}

\noindent By rearranging  the Eq. (4.37)  as
\begin{equation} \label{EQ58}
F(y'',y,\varepsilon )=y''+\pi ^{2} e^{-\varepsilon y} =Ly+N(y,\varepsilon )
\end{equation}
\noindent and using the Eqs. (2.2) and (4.3) with $\varepsilon =1$, we have
\begin{equation} \label{EQ59}
(y_{c} ^{{'} {'} } )_{n} =\pi ^{2} y_{n} -\left(y_{n} ^{{'} {'} } +\pi ^{2} \right).
\end{equation}
\noindent Without going into details here, we just give the successive iterations:

\noindent
\begin{equation}  \label{EQ60}
y_{0} =0
\end{equation}
\begin{equation}  \label{EQ61}
y_{1} =\frac{\pi ^{2} C_{0} }{2} \left(x-x^{2} \right)
\end{equation}
\begin{equation}  \label{EQ62}
y_{2} =\frac{\pi ^{2} C_{0} }{2} \left(x-x^{2} \right)-\frac{(C_{0} +C_{1} )}{24} (-x+x^{2} )\pi ^{2} \left[-12+\left(12+(-1-x+x^{2} )\pi ^{2} \right)C_{0} \right]
\end{equation}
\begin{equation} \label{EQ63}
\begin{array}{l} {y_{3} = -\frac{\pi ^{2} C_{0} }{2} (-x+x^{2} )-\frac{x\pi ^{2} }{24} (C_{0} +C_{1} )\left(-1+x\right)\left(12+C_{0} \left(-12+(-1+(-1+x)x)\pi ^{2} \right)\right)} \\ {+\frac{(C_{0} +C_{1} +C_{2} )\pi ^{2} }{720} \left[\begin{array}{l} {-360(-1+C_{0} )(-1+C_{0} +C_{1} )(-1+x)x} \\ {+30\left(-C_{1} +2C_{0} (-1+C_{0} +C_{1} )\right)(x-2x^{3} +x^{4} )\pi ^{2} } \\ {-C_{0} (C_{0} +C_{1} )x(-3+5x^{2} -3x^{4} +x^{5} )\pi ^{4} } \end{array}\right]} \end{array}
\end{equation}

\noindent Proceeding as earlier we find constants $C_{0} ,C_{1} $  and $C_{2} $:
\begin{equation} \label{EQ64}
C_{0} =0.00839960142,C_{1} =0.08178563321,C_{2} =-0.000193602314
\end{equation}
\noindent Inserting the constants into the Eq. (4.45), we obtain the approximate solution of the third order:
\begin{equation}\label{EQ65}
\begin{array}{l} {y_{3} (x)=3.134717936805843x-4.811906503098512x^{2} +4.266200757140372x^{3} } \\ {{\rm \; \; \; \; \; \; \; \; \; \; }-4.407364172185969x^{4} +2.7222682887263185x^{5} -0.9036537597026911x^{6}} \end{array}
\end{equation}

\newpage
\textit{\underbar{OPIA (1,2) }}\underbar{}

\noindent After making the relevant calculations, the algorithm
\begin{equation} \label{EQ66}
N+N_{y} (y_{c} )_{n} \varepsilon +N_{\varepsilon } \varepsilon +N_{\varepsilon y} (y_{c} )_{n} \varepsilon ^{2} +\frac{N_{\varepsilon \varepsilon } \varepsilon ^{2} }{2} +\frac{N_{yy} \varepsilon ^{2} (y_{c} )_{n}^{2} }{2} =0
\end{equation}
reduces to
\begin{equation} \label{EQ67}
(y_{c} ^{{'} {'} } )_{n} -\pi ^{2} (y_{c} )_{n} =\pi ^{2} y_{n} -\frac{\pi ^{2} }{2} y_{n}^{2} -y_{n} ^{{'} {'} } -\pi ^{2}  .
\end{equation}
Using the Eqs. (3.2), (4.49) and the initial conditions, we get
\begin{equation} \label{EQ68}
y_{0} =0
\end{equation}
\begin{equation} \label{EQ69}
y_{1} =C_{0} \left(-1+\cosh[\pi x]+\sinh[\pi x]\tanh[\frac{\pi }{2} ]\right)
\end{equation}
\begin{equation} \label{EQ70}
\begin{array}{l} {y_{2} (x)=C_{0} \left(-1+\cosh[\pi x]+\sinh[\pi x]\tanh[\frac{\pi }{2} ]\right)+\frac{{\rm e}^{-2\pi x} (-1+{\rm e}^{\pi } )}{12(1+{\rm e}^{\pi } )^{2} } (C_{0} +C_{1} )\left(-1+\coth[\pi ]\right)\times } \\ {\left[\begin{array}{l} {-6{\rm e}^{\pi x} (1+{\rm e}^{\pi } )^{2} (-1+{\rm e}^{\pi x} )(-{\rm e}^{\pi } +{\rm e}^{\pi x} )+6A{\rm e}^{\pi x} (1+{\rm e}^{\pi } )^{2} (-1+{\rm e}^{\pi x} )(-{\rm e}^{\pi } +{\rm e}^{\pi x} )} \\ {-C_{0} ^{2} \left\{\begin{array}{l} {{\rm e}^{2\pi } +{\rm e}^{3\pi } -3{\rm e}^{2\pi x} +{\rm e}^{4\pi x} -15{\rm e}^{2\pi (1+x)} -3{\rm e}^{\pi (3+2x)} -15{\rm e}^{\pi +2\pi x} +{\rm e}^{\pi +4\pi x} } \\ {+{\rm e}^{\pi +\pi x} (2+3\pi (-1+x))+{\rm e}^{3\pi x} (2-3\pi x)+3{\rm e}^{\pi +3\pi x} (4+\pi -2\pi x)+} \\ {{\rm e}^{\pi (3+x)} (2+3\pi x)+3{\rm e}^{\pi (2+x)} \left(4+\pi (-1+2x)\right)+{\rm e}^{\pi (2+3x)} \left(2-3\pi (-1+x)\right)} \end{array}\right\}} \end{array}\right]} \end{array}
\end{equation}
for OPIA(1,2). Using the Eq.(3.11) , the following values of  $C_{0} $ and $C_{1} $ are obtained:
\begin{equation} \label{EQ71}
C_{0} =-1.0286083214317654,C_{1} =2.029583070812236
\end{equation}
By using the above values, the approximate solution of the second order is:
\begin{equation} \label{EQ72}
\begin{array}{l} {y_{2} (x)=1.56204116701300+(-1.4080552742209-1.4777172978873x)\cosh[\pi x]} \\ {{\rm \; \; \; \; \; \; \; \; \; \; \; }+0.15341644132458\sinh[2\pi x]-0.15399004069623\cosh[2\pi x]} \\ {{\rm \; \; \; \; \; \; \; \; \; \; \; }+(1.163463843313+1.6111742346200355x)\sinh[\pi x]{\rm \; \; \; \; }} \end{array} .
\end{equation}

\begin{table}\caption{Comparison of absolute errors of  Example 3 at different orders of approximations.}
\begin{tabular}{|p{0.15in}|p{1in}|p{1in}|p{1in}|p{1in}|p{1in}|p{0.9in}|} \hline
       \textit{x} & \multicolumn{3}{|p{3in}|}{Absolute  errors for  OPIA(1,1)   solutions} & \multicolumn{2}{|p{2.3in}|}{Absolute  errors for  OPIA(1,2) solutions } & Exact                                solution \\ \hline
 &    $\left|y-y_{1} \right|$ &         $\left|y-y_{2} \right|$ &      $\left|y-y_{3} \right|$ &      $\left|y-y_{1} \right|$ &     $\left|y-y_{2} \right|$ &   \\ \hline
0.1 & 0.0000752784 & 0.0000608251 & 0.0000719575 & 9.05621$\times$10${}^{-7}$ & 8.8864$\times$10${}^{-7}$ & $0.269276469$ \\
0.2 & 0.0004108547 & 0.0001009657 & 0.0000183205 & 4.03657$\times$10${}^{-7}$ & 1.88504$\times$10${}^{-7}$ & $0.462340122$ \\
0.3 & 0.0000296314 & 0.0000723684 & 4.31405$\times$10${}^{-6}$ & 3.99521$\times$10${}^{-7}$ & 1.23351$\times$10${}^{-7}$ & $0.592783600$ \\
0.4 & 0.0000955682 & 0.0000135841 & 5.90773$\times$10${}^{-6}$ & 2.60399$\times$10${}^{-6}$ & 1.18017$\times$10${}^{-7}$ & $0.668371029$ \\
0.5 & 0.0002856413 & 2.90365$\times$10${}^{-5}$ & 1.28998$\times$10${}^{-6}$ & 3.05668$\times$10${}^{-6}$ & 2.73484$\times$10${}^{-7}$ & $0.693147180$ \\
0.6 & 0.0000213685 & 8.10269$\times$10${}^{-6}$ & 1.07103$\times$10${}^{-6}$ & 8.70569$\times$10${}^{-7}$ & 2.68334$\times$10${}^{-7}$ & $0.668371029$ \\
0.7 & 0.0000723646 & 9.30855$\times$10${}^{-6}$ & 1.1606$\times$10${}^{-6}$ & 5.19005$\times$10${}^{-7}$ & 9.46642$\times$10${}^{-8}$ & $0.592783600$ \\
0.8 & 0.0000108799 & 9.99237$\times$10${}^{-6}$ & 2.05027$\times$10${}^{-6}$ & 8.05111$\times$10${}^{-7}$ & 4.36073$\times$10${}^{-7}$ & $0.462340122$ \\
0.9 & 0.0005692033 & 0.000111947 & 0.0000523313 & 1.22014$\times$10${}^{-6}$ & 5.3818$\times$10${}^{-7}$ & $0.269276469$ \\ \hline
\end{tabular}
\end{table}

\begin{figure}[h]
\centerline{\includegraphics[width=0.7\textwidth]{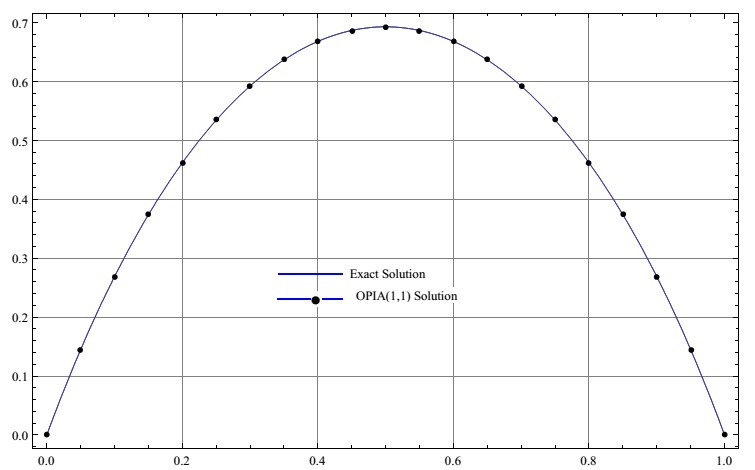}}
\caption{Comparison between the three-term OPIA(1,1) approximate solution and the exact solution for Example 3.}
\end{figure}

\noindent One  can easily observe from Table 3 and Figure 3 that the results agree very well with the exact solution.

\section{Conclusions}
In this paper, a new technique OPIM is employed for the first time to obtain a new analytic approximate solution of Bratu-type differential equations. This new method provides us with an easy way to optimally control and adjust the convergence solution series. OPIM gives a very good approximation even in a few terms which converges to the exact solution. This fact is obvious from the use of the auxiliary function $S_{n} (\varepsilon )$ which depends on \textit{n} coefficients $C_{0} ,C_{1} ,\ldots ,C_{n}.$   The results obtained in this paper confirm that the OPIM is a powerful and efficient technique for finding nearly exact solutions for differential equations which have great significance in many different fields of science and engineering.

\bibliography{BibfileBratu}{}

\end{document}